\documentclass[a4paper,11pt]{amsart}
\usepackage{hyperref}
\usepackage[cyr]{aeguill}
\newtheorem{thm}{Theorem}[section]
\newtheorem{cor}[thm]{Corollary}
\newtheorem{lem}[thm]{Lemma}
\newtheorem{prop}[thm]{Proposition}
\theoremstyle{definition}
\newtheorem{defn}[thm]{Definition}
\theoremstyle{remark}
\newtheorem{rem}[thm]{Remark}
\numberwithin{equation}{section}

\newcommand{\set}[1]{\left\{#1\right\}}
\newcommand{\clo}[1]{Cl(#1)}
\newcommand{\rclo}[2]{Cl_{#1}(#2)}
\newcommand{\bd}[1]{Bd(#1)}
\newcommand{\rbd}[2]{Bd_{#1}(#2)}

\newcommand{\Real}{\mathbb R}

\hypersetup{
    pdftitle={Recurrent Surface Homeomorphisms}, 
    pdfauthor={Boris Kolev and Marie-Christine P\'{e}rou\`{e}me}, 
    pdfsubject={MSC : 37B20, 37E30, 54H20, 57N05}, 
    pdfkeywords={recurrent, surfaces homeomorphisms, prime ends}, 
    pdflang=en, 
    pdfdisplaydoctitle=true,
    colorlinks=true,
    urlcolor=blue,
    citecolor=blue,
    linkcolor=blue,
    pdfstartview=FitH,
    pdfpagemode=None,
    }
\begin{document}

\title{Recurrent surface homeomorphisms}
\date{August 12, 1996}

\author{Boris Kolev}%
\address{CMI, 39 rue F. Joliot-Curie, 13453 Marseille cedex 13, France} %
\email{kolev@cmi.univ-mrs.fr}%

\author{Marie-Christine P\'{e}rou\`{e}me}

\thanks{The authors express their gratitude to Alexis Marin for several remarks that helped to improve this paper.}%


\begin{abstract}
An orientation-preserving recurrent homeomorphism of the
two-sphere which is not the identity is shown to admit exactly two
fixed points. A recurrent homeomorphism of a compact surface with
negative Euler characteristic is periodic.
\end{abstract}
\maketitle

\section{Introduction}

A homeomorphism $f$ of a compact metric space $(X, d)$ is
\emph{recurrent} if it admits iterates arbitrarily close to the
identity, i.e. if there exists a sequence $n_k \to + \infty$  such
that
\begin{equation*}
d\left( {f^{n_k },Id} \right) \to 0\quad \text{as} \quad k \to +
\infty
\end{equation*}

This notion is in fact independent of the metric $d$ which defines
the topology on $X$, and is invariant under topological conjugacy.
The set of recurrent homeomorphisms of $X$ includes \emph{periodic}
homeomorphisms and more generally \emph{regular} homeomorphisms
\cite{BK98,Got58}, i.e. homeomorphisms such that the sequence of
all iterates forms an equicontinuous family:

\begin{equation*}
\forall \varepsilon > 0,\;\exists \eta > 0;\quad d(x,y) < \eta
\Rightarrow d(f^{k}(x),f^{k}(y)) < \varepsilon ,\quad \forall k \in
\mathbb{Z}.
\end{equation*}

In the case of the two-sphere $\mathbb{S}^{2}$, regular
homeomorphisms have been completely classified by Ker\'{e}kj\'{a}rt\'{o}. They
are conjugate to the restriction of Euclidean isometries of the
ambient $3$-space \cite{BK98,CK94,Ker19,Ker34}.

In \cite{Ker35b}, Ker\'{e}kj\'{a}rt\'{o} asked if a non-periodic,
orientation-preserving, recurrent homeomorphism of the sphere is
always conjugate to an irrational rotation. The answer is now known
to be false. In \cite{FO92}, Fokkink and Oversteegen gave an example
of a periodic-point free, recurrent homeomorphism of the annulus
which is not conjugate to an irrational rotation. A similar
construction gives an orientation-preserving, recurrent
homeomorphism of $\mathbb{S}^{2}$ with just two fixed points (and no
other periodic points) that is not conjugate to an irrational
rotation. In that example, the homeomorphism leaves invariant, a
non-locally connected continuum $H$ with a dense orbit even if it is
not minimal: the orbit closure of each point contains a fixed point.
Besides the fact that recurrence does not imply regularity, this
example shows that the collection of orbit closures may not even
form a partition of $\mathbb{S}^{2}$. However, we will establish
here the following result:

\begin{thm}\label{thm:Main}
A non-trivial, orientation-preserving and recurrent homeomorphism of
the sphere $\mathbb{S}^{2}$ has exactly two fixed points.
\end{thm}

An immediate corollary of Theorem \ref{thm:Main} is that an
orientation-preserving, recurrent homeomorphism of the closed disc
$D^2$ which is the identity on the boundary is the identity on the
whole disc. This is the main ingredient used by Oversteegen and
Tymchatyn in \cite{OT90} to show that a recurrent homeomorphism of
the Euclidean plane $\Real^2$ is periodic. However, we do not think
that Theorem \ref{thm:Main} can be proved by a simple rephrasing of
the arguments in \cite{OT90}: one could expect to deduce Theorem
\ref{thm:Main} from the result of Oversteegen and Tymchatyn applying
the arguments of section 4 to the non-compact surface $M$ obtained
by removing the fixed points of a recurrent homeomorphism of the
sphere. The homeomorphism obtained in such a way is clearly
recurrent for the induced, restricted metric on $M$ but not
necessarily for a hyperbolic metric on $M$. Notice that the
analogous statement is false in a simpler case: if one removes one
or two fixed points of an irrational rotation of the sphere, the
induced homeomorphism on $M$ (the plane or the open annulus) is
certainly not recurrent for the euclidean metric, otherwise it would
be periodic by \cite{OT90}. Moreover, our proof relies on a simple
property: the fact that a homeomorphism is recurrent is inherited by
the induced homeomorphism on prime ends. Up to our knowledge, this
elementary idea has not yet been pointed out in the literature.

Even if he does not seem to have published any proof of it,
Ker\'{e}kj\'{a}rt\'{o} announced in \cite{Ker35b} that a recurrent homeomorphism
of a closed, orientable surface of genus $g > 1$ is periodic. This
result is a corollary of Theorem \ref{thm:Main} as we shall see in
section 4. As noted by Ker\'{e}kj\'{a}rt\'{o} in \cite{Ker35a}, this implies in
particular that a homeomorphism of a compact surface with negative
Euler characteristic cannot be topologically transitive provided it
has at least one regular point.

The paper is organized as follows: in section 2, we review required
material from \emph{prime ends'} theory. For the proofs in this
section, we refer the reader to the works of Epstein \cite{Eps81} or
Mather \cite{Mat82}. The proof of Theorem \ref{thm:Main} is given in
section 3. Section 4 contains several corollaries concerning
recurrent homeomorphisms of other surfaces orientable or not.


\section{Prime ends}

Let $U$ be an open, connected subset of the sphere $\mathbb{S}^{2}$.
For a set $A \subset U$ we let $\clo{A}$ be its closure in
$\mathbb{S}^{2}$, $\rclo{U}{A}$ be its closure in $U$, $\bd{A}$ be
its boundary in $\mathbb{S}^{2}$ and $\rbd{U}{A}$ be its boundary in
$U$.

\begin{defn}
A \emph{topological chain} is a sequence $\Omega = (\Omega_{i} )_{i
\in {\mathbb{N}}}$ of open subsets of $U$ such that:
\begin{enumerate}
    \item $\Omega_{i}$ and $\rbd{U}{\Omega_{i}}$ are non-empty connected sets;
    \item $\rclo{U}{\Omega_{i + 1}} \subset \Omega_{i}$
    \item $\clo{{\rbd{U}{\Omega_{i}} } } \cap \clo{
{\rbd{U}{\Omega_{j}} } } = \emptyset \quad if\quad i \ne j$;
    \item there exists a point $x \in \clo{U}$, the \emph{principal
point} of $\Omega$, such that $x \notin \clo{ \rbd{U}{\Omega_{i}}
}$, for all $i$ and $\rbd{U}{\Omega_{i}} \to x$ as $i \to \infty$.
\end{enumerate}
\end{defn}

Two chains $\Omega $ and $\Omega^{\prime}$ are \emph{equivalent},
denoted by $\Omega \simeq \Omega^{\prime}$, if they divide each
other, i.e. if for all $i$ there exists $j$ such that
$\Omega^{\prime}_j \subset \Omega_{i} $ and $\Omega _j \subset
\Omega^{\prime}_i $. The set of all topological chains equivalent to
a given chain $\Omega $ is called the \emph{prime point} defined by
$\Omega$, and we write $e= [\Omega ]$, or $\Omega \in e$. We denote
the set of all the prime points of $U$ by  $\hat{U}$.

The set of \emph{principal points} of chains $\Omega \in e$ is
called the \emph{principal set} of $e$, which we denote by $X(e)$.
Let $e$ be a prime point, and let $\Omega $ be any topological chain
that defines $e$. The set
\begin{equation*}
    Y(e) =  \bigcap_{i \in {\mathbb{N}}} \clo{\Omega_{i}}
\end{equation*}
depends only on $e$, and is called the \emph{impression} of $e$. It
can be shown \cite{Eps81,Mat82} that $Y(e)$ and $X(e)$ are continua
(i.e. non empty, compact, connected sets).

Given $x \in U$, choose any sequence of open discs $(D_{i})_{i \in
{\mathbb{N}}}$ such that
\begin{equation*}
    \clo{D_{i + 1}} \subset D_{i} \subset U \quad \text{and} \quad \bigcap_{i \in {\mathbb{N}}} \clo{D_{i}} = \set{x}.
\end{equation*}
The sequence $(D_{i})_{i \in {\mathbb{N}}}$ is a topological chain,
and the prime point $\omega(x)$ which it defines depends only on
$x$, and not of the particular decreasing sequence $(D_{i})_{i \in
{\mathbb{N}}}$. Moreover, we have that
\begin{equation*}
    X(\omega(x)) = Y(\omega(x)) = \set{x} .
\end{equation*}
This allows us to define a one-to-one map $x \mapsto \omega(x)$ from
$U$ into $\hat{U}$, and to consider $U$ as a subset of $\hat{U}$. A
prime point which is not in $\omega(U)$ is called a \emph{prime
end}, and it is characterized by the fact that $Y(e) \subset \bd{U}$
(see Lemma 2.5 of \cite{Eps81}).

A point $x \in \bd{U}$ is \emph{accessible} (from $U$) provided that
there exists a path $\gamma :[0,1[ \to U$ such that $\gamma \left( t
\right) \to x$ as $t \to 1$. The set of accessible points is dense
in $\bd{U}$. If $\mathbb{S}^{2}\setminus U$ is the union of a finite
number of non-degenerate continua, we have the following lemma.

\begin{lem}\label{lem:accessible}
For each accessible point $x \in \bd{U}$, there exists a prime end
$e \in \hat {U}\setminus U$ such that $X(e) = \set{x}$
\end{lem}

For any open set $V \subset U$, we define
\begin{equation*}
\hat {V} = \set{[\Omega] \in \hat{U}; \quad \Omega_{i} \subset
V,\quad \text{for some} \; i}.
\end{equation*}
With these notations, we get
\begin{equation*}
    \hat {V} \cap \hat {W} = \widehat{V \cap W},
\end{equation*}
and
\begin{equation*}
    \hat {V} \cup \hat {W} \subset \widehat{V \cup W}.
\end{equation*}
Therefore, the family $(\hat {V})_{V \subset U}$ is a basis for a
topology $\tau $ on $\hat{U}$. It follows from the definition of
this topology that:

\begin{enumerate}
    \item if $e = [\Omega]$ then $(\hat{\Omega}_{i})_{i\in\mathbb{N}}$ is
a basis of neighborhoods of $e$,
    \item $\omega :U \to \hat{U}$ is a continuous, open map, and
    $\omega(U)$ is a dense open subset of $\hat {U}$.
\end{enumerate}

The following theorem is a generalization of a theorem due
originally to Carath\'{e}odory~\cite{Car13} for a simply connected
domain.

\begin{thm}[\cite{Eps81}, Theorem 6.6]\label{thm:Epstein}
Let $U$ be an open subset of the sphere $\mathbb{S}^{2}$ such that
$\mathbb{S}^{2}\setminus U$ is the union of a finite number $p$ of
non-degenerate continua. Then, $\hat {U}$, the \emph{prime end
compactification} of $U$ is homeomorphic to a compact surface of
genus $0$ with $p$ boundary components.
\end{thm}

Suppose now that we are given a self-homeomorphism $f$ of the pair
$\mathbb{S}^{2}$, such that $f(U) = U$. For each topological chain
$\Omega$ with principal point $x$, $f(\Omega) = \left( f(\Omega_{i})
\right)_{i \in {\mathbb{N}}}$ is a topological chain with principal
point $f(x)$ and $\Omega^{\prime} \simeq \Omega $ if and only if
$f(\Omega^{\prime}) \simeq f(\Omega)$. Hence $f$ induces a
homeomorphism $\hat{f}$ on $\hat{U}$ such that:

\begin{enumerate}
    \item for each $V \subset U$, $\hat{f}(\hat{V}) = \widehat{f(V)}$,
    \item for each $e \in \hat{U}$, $X(\hat{f}(e)) = f(X(e))$,
    \item for each $e \in \hat{U}$, $Y(\hat{f}(e)) = f(Y(e))$.
\end{enumerate}

Clearly, if $x \in U$ is a fixed point of $f$ then $\omega \left( x
\right)$ is a fixed point for $f$. But in general, nothing can be
said about fixed prime ends from the existence of fixed points of
$f$ lying on $\bd{U}$. Moreover, there are simple examples which
show that $f$ may have fixed points on $\bd{U}$ although $\hat{f}$
has no fixed point at all. However, if $e$ is a fixed prime end and
$X(e) = \set{x}$ then $f(x) = x$ since
\begin{equation*}
\set{f(x)} = X(\hat{f}(e)) = X(e) = \set{x} .
\end{equation*}

The following result, which we state without proof, is a consequence
of Lemma~\ref{lem:accessible} and of the density of accessible
points in $\bd{U}$.

\begin{cor}\label{cor:boundary}
Let $f$ be a homeomorphism of $\mathbb{S}^{2}$ and $U$ an invariant
open set such that $\mathbb{S}^{2}\setminus U$ is the union of a
finite number of non-degenerate continua. If $\hat {f} = \hat{Id}$
on $\hat{U}\setminus \omega(U)$, then  $f = Id$ on $\bd{U}$.
\end{cor}

The following lemma is new. The fact that a homeomorphism is
recurrent is inherited by the induced homeomorphism on prime ends
seems not to have been pointed out in the literature, before.

\begin{lem}\label{lem:main}
Let $f$ be a recurrent homeomorphism of the sphere $\mathbb{S}^{2}$
and $U$ an invariant open set, such that $\mathbb{S}^{2}\setminus U$
is the union of a finite number of non-degenerate continua. Then,
the homeomorphism $\hat{f}$ induced by $f$ on $\hat {U}$ is
recurrent.
\end{lem}

Before giving the proof of this result, we recall an elementary
lemma of topology widely used in the theory of prime ends
\cite{Mat82}.

\begin{lem}\label{lem:topology}
Let $A$ and $B$ be open subsets of a topological space $X$. Suppose
that $X$, $A$, $B$, $\rbd{X}{A}$ and $\rbd{X}{B}$ are non-empty
connected sets, and that $\rbd{X}{A} \cap \rbd{X}{B} = \emptyset$.
Then exactly one of the following situations holds:
\begin{enumerate}
    \item $X = A \cup B$,
    \item $\rclo{X}{A} \cap \rclo{X}{B} = \emptyset$,
    \item $\rclo{X}{A} \subset B$,
    \item $\rclo{X}{B} \subset A$.
\end{enumerate}
\end{lem}

\begin{proof}[Proof of Lemma~\ref{lem:main}] We fix a metric on the
compact surface of finite type $\hat {U}$, and to avoid any
confusion with the metric $d$ on $U$, we denote it by $\alpha$.
Choose a sequence $(n_k)$ such that
\begin{equation*}
    d(f^{n_k },Id) \to 0
\end{equation*}
on $\clo{U}$. Let $\varepsilon
> 0$, $e \in \hat {U}$ and $B_\alpha \left( {e,\varepsilon }
\right)$ be the ball of radius $\varepsilon$ around $e$. For any
chain $\Omega = (\Omega_{i})_{i\in\mathbb{N}}$ that defines $e$, we
have $\hat {\Omega }_i \subset B_\alpha \left( {e,\varepsilon }
\right)$, for large enough $i$. Hence we can assume that $\hat
{\Omega }_1 \subset B_\alpha \left( {e,\varepsilon } \right)$.

For $k$ large enough, say for $k$ greater than some integer
$k(e,\varepsilon)$, the hypothesis of Lemma~\ref{lem:topology}
applies to the sets $X = U$, $A = \Omega _1$ and $B = f^{n_k }\left(
{\Omega _2 } \right)$. The only possibility among the four
alternative is the fourth one and so
\begin{equation*}
    f^{n_k }(\Omega_{2}) \subset \Omega_{1} \quad \text{for all}
    \quad k \geq k(e,\varepsilon).
\end{equation*}
Therefore, we have found a neighbourhood $\mathcal{V}_{e} =
\widehat{\Omega}_{2}$ of $e$, and an integer $k(e,\varepsilon)$ so
that
\begin{equation*}
    \mathcal{V}_{e}\subset B_{\alpha}(e,\varepsilon) \quad
    \text{and} \quad \hat{f}^{n_{k}}(\mathcal{V}_{e})\subset
    B_{\alpha}(e,\varepsilon),
\end{equation*}
for all $k \geq k(e,\varepsilon)$. Since $\hat {U}$ is compact, only
a finite number of $\mathcal{V}_e$, say $\mathcal{V}_{e_{1}} ,
\mathcal{V}_{e_{2}} , \cdots , \mathcal{V}_{e_{r}}$, are necessary
to cover $\hat{U}$, and hence, for
\begin{equation*}
    k \geq \max \set{k(e_{1},\varepsilon), k(e_{2},\varepsilon), \dotsc ,
    k(e_{r},\varepsilon)},
\end{equation*}
we have $\alpha \left( {f^{n_k },Id} \right) \le 2\varepsilon$. In
other words, $\hat {f}$ is recurrent.
\end{proof}

\section{Proof of the main proposition}

The proof of Theorem~\ref{thm:Main} relies heavily on the
construction, for every fixed point of $f$, of arbitrarily small,
invariant, non degenerate continua containing the fixed point: the
so-called \emph{Birkhoff construction}.

\begin{lem}[Birkhoff Construction]\label{lem:birkhoff}
Let $f$ be a \emph{recurrent homeomorphism} of the sphere
$\mathbb{S}^{2}$, $x_{0}$ a fixed point of $f$ and $D$, a
topological disc, containing $x_{0}$, and bounded by a simple closed
curve $c$. There exists a simply-connected continuum $K$, such that:
\begin{enumerate}
    \item $x_{0} \in K \subset \clo{D}$,
    \item $f\left( K \right) = K$,
    \item $K \cap c \ne \emptyset$.
\end{enumerate}
\end{lem}

\begin{proof}
We define inductively a decreasing sequence of Jordan domains,
$D_{n}$, by the following properties
\begin{enumerate}
    \item $D_{0} = D$,
    \item $D_{n + 1}$ is the connected component of $x_{0}$ in $f(D_{n}) \cap D_{0}$.
\end{enumerate}

That each $D_n$ is a Jordan domain can be established by induction
on $n$: if  $f(D_{n}) \subset D_{0}$ or $D_{0} \subset f(D_{n})$,
there is nothing to prove. Otherwise, the two simple closed curves
$\bd{D_{0}}$ and $f(\bd{D_{n}})$ have at least two common points and
$D_{n + 1}$ is a Jordan domain (cf. Theorem 16.3 of \cite{New92}).

As the intersection of decreasing discs, the set
\begin{equation*}
    K = \bigcap_{n \in {\mathbb{N}}} \clo{D_{n}}
\end{equation*}
is a simply-connected continuum. It is contained in $\clo{D}$ and
contains $x_{0}$. Moreover $K \subset f(K)$ and \emph{since $f$ is
recurrent} this implies $f(K) = K$.

It remains to show that $K \cap c \ne \emptyset$, the proof of which
is a consequence of the \emph{intersecting curve property}, verified
by a recurrent homeomorphism. Suppose on the contrary that $K \cap c
= \emptyset$. In that case, we can find an integer $n_{0} \ge 0$
such that
\begin{equation*}
    \clo{D_{n_{0} + 1}} \cap c = \emptyset
\end{equation*}
and hence
\begin{equation*}
    \bd{D_{n_{0} + 1}} \subset \bd{f(D_{n_{0}})},
\end{equation*}
which leads to
\begin{equation*}
    D_{n_{0} + 1} = f(D_{n_{0}}).
\end{equation*}
But then, we have
\begin{equation*}
    \bd{D_{n_{0} + 2}} \subset \bd{f(D_{n_{0} + 1})} = \bd{f^{2}(D_{n_{0}})},
\end{equation*}
and hence
\begin{equation*}
    D_{n_{0} + 2} = f^{2}(D_{n_{0}}).
\end{equation*}
Iterating the process inductively, we obtain that
\begin{equation*}
    D_{n_{0} + i} = f^{i}(D_{n_{0}}),
\end{equation*}
for all $i \ge 0$, and hence that
\begin{equation*}
    K = \bigcap_{i\in\mathbb{N}} \clo{f^{i}(D_{n_{0}})} .
\end{equation*}
But then $K \cap c_{n_{0} } \ne \emptyset$, otherwise, we could find
$k \in \mathbb{N}$ such that
\begin{equation*}
    \clo{f^{k}(D_{n_{0}})} \cap c_{n_{0}} = \emptyset ,
\end{equation*}
which is not possible since $f$ is recurrent. Moreover, since
\begin{equation*}
c_{n_{0}} \subset c \cup f(c) \cup \dotsb \cup f^{n_{0}}(c) ,
\end{equation*}
and $f(K) = K$, the existence of a point in $K \cap c_{n_{0}}$ would
lead to the existence of a point in $K \cap c$, which gives a
contradiction and completes the proof.
\end{proof}

\begin{prop}\label{prop:connected}
Let $f$ be a recurrent, orientation-preserving homeomorphism of the
sphere $\mathbb{S}^{2}$. If $f$ has three fixed point, then $Fix(f)$
is connected.
\end{prop}

\begin{proof}
Suppose on the contrary that $Fix(f)$ is not connected. Then, we can
write
\begin{equation*}
    Fix(f) = X \cup Y ,
\end{equation*}
where $X$ and $Y$ are two disjoint, non-empty closed sets. Let $x
\in X$ and $y \in Y$. We can find a Jordan curve $c$, in the
complement of $Fix(f) = X \cup Y$, which separates $x$ from $y$ (see
\cite{New92} or \cite{Why64}). Let $D_{x}$ (resp. $D_{y}$) be the
component of $\mathbb{S}^{2}\setminus c$ which contains $x$ (resp.
$y$). By hypothesis, $f$ has a third fixed point, and hence, at
least one component of $\mathbb{S}^{2}\setminus c$, say $D_{x}$,
contains two fixed points $x$ and $z$. Let $J \subset D_{x}$ be a
simple closed curve which separates $x$ and $z$ in $D_{x}$, and let
$\Delta_{x}$ be the component of $\mathbb{S}^{2}\setminus J$ that
contains $x$. Using the construction of Lemma~\ref{lem:birkhoff},
first for $x$ and $\Delta_{x}$ and then for $y$ and $D_{y}$, we
obtain two disjoint, invariant, non-degenerate, simply connected
continua $K_{x}$ and $K_{y}$ such that
\begin{equation*}
    \emptyset \ne K_{y} \cap c \subset \bd{K_{y}} \cap c \subset \mathbb{S}^{2}\setminus Fix(f)
\end{equation*}

The set $U = \mathbb{S}^{2}\setminus (K_{x}\cup K_{y})$ is an open
topological annulus which contains $z$ and which is invariant under
$f$. According to Theorem~\ref{thm:Epstein}, $\hat{U}$ is
homeomorphic to a closed annulus and due to Lemma~\ref{lem:main},
the induce homeomorphism $\hat{f}$ is recurrent. Moreover, $\hat{f}$
is orientation-preserving, boundary-preserving and has a fixed point
$\omega(z)$.

In general, a lift $F$ of $\hat{f}$ to $\Real\times [0,1]$ needs not
be recurrent. However if we choose a lift $F$ which has a fixed
point, it is recurrent. Therefore, on each line $\Real\times
\set{i}$, ($i = 0,1$), $F$ is the identity. Hence $\hat{f}$ itself
is the identity on $\hat{U}\setminus U$ and according to
Corollary~\ref{cor:boundary}, $f$ is the identity on $\bd{U}$. But
$\bd{U}$ meets $c$ which lies in $\mathbb{S}^{2}\setminus Fix(f)$.
This gives a contradiction and completes the proof.
\end{proof}

\begin{proof}[Proof of Theorem~\ref{thm:Main}]
Let $f$ be an orientation-preserving homeomorphism of
$\mathbb{S}^{2}$. According to Brouwer's fixed point theorem for the
sphere, $f$ has at least one fixed point $x$. Moreover, if $f$ is
recurrent then $f$ induces an orientation-preserving homeomorphism
of
\begin{equation*}
    \mathbb{S}^{2}\setminus \set{x}\simeq \Real^2
\end{equation*}
for which each point is recurrent. According to Brouwer's lemma on
translation arcs \cite{Bro84,Gui94}, which asserts that all points
of an orientation-preserving, fixed point free homeomorphism of the
plane are wandering, $f$ has necessarily a second fixed point $y$.

Suppose that $f$ has a third fixed point. According to
Proposition~\ref{prop:connected}, $Fix(f)$ is connected. Hence, each
component of $\mathbb{S}^{2}\setminus Fix(f)$ is a simply connected
domain of $\mathbb{S}^{2}$ homeomorphic to the plane $\Real^{2}$.

According to a result of Brown and Kister \cite{BK84}, an
orientation-preserving homeomorphism of the sphere leaves invariant
each component of the complement of its fixed point set.

Hence, if $Fix(f) \ne \mathbb{S}^{2}$ there exists a component $U
\simeq \Real^2$ of $\mathbb{S}^{2}\setminus Fix(f)$ such that $f(U)
= U$ and $Fix(f_{U}) = \emptyset$. But each point of $U$ is
recurrent under $f$, this is once again a contradiction with
Brouwer's lemma on translation arcs. This completes the proof of
Theorem~\ref{thm:Main}.
\end{proof}

\begin{cor}\label{cor:reversing}
A recurrent orientation-reversing homeomorphism of the sphere which
has a fixed point is an involution.
\end{cor}

\begin{proof} Let $f$ be a recurrent, orientation-reversing
homeomorphism of the sphere. Suppose that $f$ has a fixed point $x$
but $f^2 \ne Id$. According to Theorem \ref{thm:Main}, $f^2$ has
only two fixed points $x$ and $y$ and since $f(x)=x$, we must have
\begin{equation*}
    Fix(f) = \set{x,y}.
\end{equation*}
The Birkhoff's construction may be applied to obtain two distinct,
$f$-invariant, simply connected continua $K_{x}$ and $K_{y}$. As
described in sections 2 and 3, $f$ induces a recurrent,
orientation-reversing homeomorphism $\hat{f}$ on the prime end
compactification of $\mathbb{S}^{2}\setminus (K_x \cup K_y)$ which
is a closed annulus.

Since $\hat {f}$ preserves the boundary component of this annulus
and reverse the orientation, it must have a fixed point according to
the Lefschetz trace formula. Hence, $\hat {f}^2$ must be the
identity on the boundary of the annulus. But this implies that $f^2
= Id$ on $\bd{K_x \cup K_y}$ which leads to a contradiction as in
the proof of Theorem \ref{thm:Main}.
\end{proof}

\begin{cor}\label{cor:disc}
A non trivial, recurrent, orientation-preserving homeomorphism of
the disc has a unique fixed point. A recurrent,
orientation-reversing homeomorphism of the disc is an involution.
\end{cor}

\begin{proof} Taking the double of a recurrent,
orientation-preserving homeomorphism of the disc, we obtain at once
the first part of Corollary~\ref{cor:disc} as a consequence of
Theorem~\ref{thm:Main}. Let us now consider a recurrent,
orientation-reversing homeomorphism $f$ of the disc. As an
orientation-reversing homeomorphism of the circle, the restriction
of $f$ to the boundary of the disc has a fixed point and hence $f^2
= Id$ on the boundary of the disc. By previous considerations, this
leads to $f^2 = Id$ and completes the proof.
\end{proof}


\section{Recurrent homeomorphisms of surfaces}

Let $M^2$ be a closed orientable surface of genus $g > 1$ and $\pi
:\tilde {M}^2 \to M^2$ the universal cover of $M^2$. We can identify
$\tilde {M}^2$ either to the euclidean plane $\Real^2$ or to the
Poincar\'{e} disc $\mathbb{D}$ in such a way that $M^2$ is homeomorphic
to the quotient of $M^2$ by a discrete subgroup $\Gamma $ of
euclidean translations or hyperbolic isometries according to whether
$M^2$ is $\mathbb{R}^2$ or $\mathbb{D}$. The metric we shall use on
$M^2$ is the quotient metric on $M^2 / \Gamma $ defined by

\begin{equation*}
    d\left( {\pi \left( x \right),\pi \left( y \right)} \right) =
\mathop {\inf }\limits_{g,h \in \Gamma } \tilde {d}\left(
{g.x,h.y} \right),
\end{equation*}
where $\tilde{d}$ is the natural metric on $\tilde{M}^2$.

There is another metric on $\tilde{M}^2$ that we shall use in the
following, namely the spherical metric. The Alexandroff
compactification $\tilde {M}^2 \cup \left\{ \infty \right\}$, is
homeomorphic to the sphere $\mathbb{S}^{2}$. The standard metric of
$\mathbb{S}^{2}$ induces a metric $\partial$ on $\tilde {M}^2$ that
we call the spherical metric. These two metrics $\tilde{d}$ and
$\partial$ are not uniformly equivalent on $\tilde {M}^2$ but
$Id:\left( {\tilde {M}^2,\tilde {d}} \right) \to \left( {\tilde
{M}^2,\partial } \right)$ is uniformly continuous.

\begin{lem}\label{lem:surfaces}
Let $f$ be an orientation-preserving recurrent homeomorphism of a
closed orientable surface $M^2$ of genus $g > 1$. If $f$ has a fixed
point and acts trivially on $\pi_{1}(M^2)$ then $f = Id$.
\end{lem}

\begin{proof} A lift of a recurrent homeomorphism of $M^2$ needs
not be recurrent in general. However, a lift $\tilde{f}$ which has a
fixed point is recurrent for the metric $d$. The fact that
$Id:(\tilde {M}^2,\tilde {d}) \to (\tilde {M}^2,\partial)$ is
uniformly continuous is enough to ensure that the extension of
$\tilde{f}$ to $\mathbb{S}^{2}$ (letting $\tilde{f}(\infty) =
\infty$) is recurrent for the metric $\partial$. Moreover if $f$
acts trivially on $\pi_{1}(M^2)$, $\tilde{f}$ has infinitely many
fixed points since it commutes with all covering translations and
hence $\tilde{f} = Id$ according to Theorem~\ref{thm:Main}.
\end{proof}

\begin{cor}\label{cor:surfaces}
A recurrent homeomorphism of a compact surface with negative Euler
characteristic is periodic.
\end{cor}

\begin{proof} We first note that if the boundary of $M^2$ is not
empty, the natural extension of $f$ to the double $DM^2$ of $M^2$ is
still recurrent and since $\chi(DM^2) = 2\chi(M^2)$, we are reduced
to prove Corollary~\ref{cor:surfaces} for closed surfaces. Moreover,
by passing to the orientation covering of $M^2$ and considering
$f^2$ instead of $f$ if necessary, we may assume that $M^2$ is
orientable with genus $g > 1$ and that $f$ is
orientation-preserving. Recall then that every homeomorphism of
$M^2$ which is close enough to the identity must be homotopic to the
identity. Hence, since $f$ is recurrent we can find a positive
integer $n$ such that $f^n$ is homotopic to the identity. The
Lefschetz number of $g = f^n$ is thus $L(g) = \chi(M^2) < 0$ and $g$
has a fixed point. According to Lemma~\ref{lem:surfaces}, we must
have then $g = Id$ which completes the proof.
\end{proof}

\begin{rem}
We emphasize on the fact that Corollary~\ref{cor:surfaces} is false
for surfaces with non negative Euler characteristic. Indeed, the
example of \cite{FO92} can be modified to exhibit non-regular,
recurrent homeomorphisms on the projective plane, the torus, the
Klein bottle, the annulus, the M\"{o}bius strip or the disc. However,
for surfaces with $\chi(M^2) = 0$, namely the torus, the Klein
bottle, the annulus and the M\"{o}bius strip, the proof of
Corollary~\ref{cor:surfaces} shows that the existence of a periodic
point for a recurrent homeomorphism implies that this homeomorphism
is periodic itself.
\end{rem}

\begin{cor}
A non-trivial, recurrent homeomorphism of the projective plane
$\mathbb{P}^{2}$ admit a unique fixed point or is an involution.
\end{cor}

\begin{proof}
We shall use as model for the real projective plane
$\mathbb{P}^{2}$, the quotient of the sphere $\mathbb{S}^{2}$ by the
involution $\theta :x \mapsto - x$. Let $f$ be a recurrent
homeomorphism of the projective plane $\mathbb{P}^{2}$. $f$ has two
lifts on $\mathbb{S}^{2}$ which commute with $\theta$. One of them
$\tilde {f}_{+}$ is orientation preserving and the other $\tilde
{f}_{-} = \theta \circ \tilde {f}_{+}$ is orientation-reversing.
$\tilde {f}_{+}$ has a fixed point and is therefore recurrent as we
have already noticed. Therefore, if $f$ is not trivial, $\tilde
{f}_{+}$ has exactly two fixed point $x$ and $\theta(x)$ which
project down to one fixed point of $f$. Hence, the existence of a
second fixed point  of $f$ implies that $\tilde {f}_{-}$ has a fixed
point and is therefore recurrent. It must satisfy $\tilde {f}_{-}^2
= Id$ according to Corollary~\ref{cor:reversing}. This shows that
$f$ must be an involution and completes the proof.
\end{proof}

\begin{rem}
Recall that a \emph{regular} point for a homeomorphism $f$ on a
compact metric space $(X, d) $ is a point at which the family of all
the iterates $\set{f^{n}; \; n \in \mathbb{Z}}$ is equicontinuous.
As noted in the introduction, Corollary~\ref{cor:surfaces} is the
key argument to show that a homeomorphism of a compact surface with
negative Euler characteristic which has a regular point cannot be
topologically transitive and hence cannot be ergodic with respect to
an invariant measure $\mu$ such that $\mu(U) > 0$ for all open sets
$U$. Indeed, the orbit of a regular point cannot be a dense set,
otherwise the homeomorphism is recurrent, and hence periodic, which
gives a contradiction. But the presence of a regular point with a
non dense orbit implies existence of proper, non empty invariant
open sets.
\end{rem}


\end{document}